\documentclass[11pt]{article}
\usepackage{latexsym,amssymb}
\usepackage{amsmath}
\usepackage{amscd}

\input xypic
\newcommand{\rar}{\rightarrow}
\newcommand{\lar}{\longrightarrow}

\newtheorem{Theorem}{Theorem}[section]

\newtheorem{Corollary}[Theorem]{Corollary}
\newtheorem{Proposition}[Theorem]{Proposition}
\newtheorem{Conjecture}[Theorem]{Conjecture}

\newtheorem{Remark}[Theorem]{Remark}
\newtheorem{Example}[Theorem]{Example}
\newtheorem{Definition}[Theorem]{Definition}

\def\Symi{\mbox{$\mathcal S$}}
\def\demo{\noindent{\bf Proof. }}
\def\QED{\hfill$\Box$}
\def\Sym{\mbox{\rm Sym}}
\def\codim{\mbox{\rm codim }}
\def\red{\mbox{\rm r}}
\def\depth{\mbox{\rm depth }}

\def\ann{\mbox{\rm ann }}
\newcommand{\Rees}{\mbox{${\mathcal R}$}}

\def\pp{{\mathbb P}}

\def\F{{\mathcal{ F}}}

\begin{document}

\title{\LARGE\sc  On the Homology of Two-Dimensional Elimination 
\vspace{-1mm}}\footnotetext{AMS 2000 {\it Mathematics Subject
Classification}. Primary 13H10; Secondary 13F20, 05C90.}

\author{
{\normalsize\sc Jooyoun Hong}
\vspace{-0.75mm}\\
{\small Department de Mathematics}\vspace{-1.4mm} \\
{\small University of California, Riverside}\vspace{-1.4mm}\\
{\small 900 Big Springs Drive}\vspace{-1.4mm} \\
{\small Riverside, CA 92521}\vspace{-1.4mm} \\
{\small e-mail: {jyhong@math.ucr.edu}}\\
\and
{\normalsize\sc Aron Simis\thanks{Partially supported by CNPq, Brazil.}}
\vspace{-0.75mm}\\
{\small Departamento de Matem\'atica}\vspace{-1.4mm} \\
{\small Universidade Federal de Pernambuco}\vspace{-1.4mm}\\
{\small 50740-540 Recife, PE, Brazil}\vspace{-1.4mm} \\
{\small e-mail: {aron@dmat.ufpe.br}}\\
\and
{\normalsize\sc Wolmer V. Vasconcelos\thanks{Partially
 supported by the NSF.}}
\vspace{-0.75mm}\\
{\small Department of Mathematics}\vspace{-1.4mm} \\
{\small Rutgers University}\vspace{-1.4mm}\\
{\small 110 Frelinghuysen Rd}\vspace{-1.4mm}\\
{\small Piscataway, New Jersey 08854-8019}\vspace{-1.4mm} \\
{\small e-mail: {vasconce@math.rutgers.edu}}\vspace{4mm}
}

\maketitle

\begin{abstract}

\noindent We study birational maps with  empty base locus defined by
almost complete intersection ideals. Birationality is shown to be
expressed by the equality of two Chern numbers. We  provide a
relatively effective method of their calculation in terms of certain
Hilbert coefficients. In dimension two the structure of the
irreducible ideals leads naturally to the calculation of Sylvester
determinants via a computer-assisted method. For degree at most $5$
we produce the full set of defining equations of the base ideal. The
results answer affirmatively some questions raised by D. Cox
(\cite{Cox}).
\end{abstract}

\date{\empty }

\newpage

\section{Introduction}
 Let $R$ be a Noetherian ring
 and $\mathbf{f}=\{ f_1, \ldots, f_m\}$  a set of elements of $R$.
Such sets are the ingredients of rational maps between affine and
other spaces. At the cost of losing some definition, we choose to
examine them  in the setting of the ideal $I$ they generate.
Specifically, we consider
 the presentation of the Rees algebra of $I$
\[ 0 \rar M \lar S = R[T_1, \ldots, T_m] \stackrel{\varphi}{\lar}
R[It] \rar 0,  \quad T_i \mapsto f_it .\] The context of  Rees
algebra theory allows for the examination of the syzygies of the
$f_i$ but also of the  relations of all orders, which are carriers
of analytic information.

 We set $\Rees = R[It]$ for
the Rees algebra of $I$. The ideal $M$ will be referred to as the
{\em equations} of the $f_j$, or by abuse of terminology, of the
ideal $I$. If $M$ is generated by forms of degree $1$, $I$ is said
to be of linear type
 (this is independent of the set of generators). The Rees algebra
 $R[It]$ is then the symmetric algebra $\Symi=\Sym(I)$ of $I$.
Such is the case when
 the $f_i$ form a regular sequence, $M$ is then generated by the
 Koszul forms $f_iT_j-f_jT_i$, $i<j$.
We will treat mainly {\em almost complete intersections}  in a
Cohen-Macaulay ring $R$, that is, ideals of codimension $r$
generated by $r+1$ elements. Almost exclusively, $I$ will be an
ideal of finite co-length in a local ring, or in a ring of
polynomials over a field.

Our focus on $\Rees$ is shaped by the following fact. The class of
ideals $I$ to be considered will have the property that both  its
symmetric algebra $ \Symi $ and the normalization $\Rees'$ of
$\Rees$ have amenable properties, for instance, one of them (when
not both) is Cohen-Macaulay. In such case,  the diagram
\[ \Symi \twoheadrightarrow \Rees \subset \Rees'\]
gives a convenient dual platform from which to examine $\Rees$.

There are  specific motivations for looking at (and for) these
equations. In order to describe our results in some detail, let us
indicate their contexts.

\begin{itemize}

\item[{\rm (i)}] Ideals which are almost complete intersections
occur in some of the more notable birational maps and in geometric
modelling (\cite{BuJou}, \cite{BuChaJou}, \cite{BuCoxDan},
 \cite{CoxSeCh}, \cite{CoxGoZh}, \cite{Cox0}, \cite{Cox},
\cite{Dand}, \cite{SGH}, \cite{bir2003}, \cite{ram2}).

\item[{\rm (ii)}]   It is possible
 interpret  questions of birationality of
certain maps as an interaction between the Rees algebra of the ideal
and its special fiber. The mediation is carried by the first Chern
coefficient of the associated graded ring of $I$. In the case of
almost complete intersections
 the analysis is more tractable, including the construction of
suitable algorithms.

\item[{\rm (iii)}]
At a recent talk in Luminy (\cite{Cox}), D. Cox raised several
questions about the character of the equations of Rees algebras in
polynomial rings in two variables.
 They are addressed in Section~\ref{dimtwo} as
 part of a general program of devising algorithms that produce
all the equations of an ideal, or at least some distinguished
polynomial (e.g. the `elimination equation' in it) (\cite{BuJou},
\cite{Jou97}).

\end{itemize}

We now describe our results. Section~\ref{Prelim} is an assemblage
for the ideals treated here of basics on symmetric and Rees
algebras, and on their Cohen-Macaulayness. We also introduce the
general notion of a Sylvester form in terms of contents and
coefficients in a polynomial ring over a base ring. This is
concretely taken up in Section~\ref{dimtwo} when the base ring is a
polynomial ring in $2$ variables over a field.

In Section~\ref{Ratparam} we examine  the connection between typical
algebraic invariants and the geometric background of rational maps
and their images. Here, besides the dimension and the degree of the
related algebras, we also consider the Chern number $e_1(I)$ of an
ideal. In particular we explain a criterion for a rational map to be
birational in terms of an equality of two such Chern numbers,
provided the base locus of the map is empty and defined by an almost
complete intersection ideal.

In Section~\ref{dimtwo}, we discuss the role of irreducible ideals
in producing Sylvester forms. Of a general nature, we describe a
method to obtain an irreducible decomposition of ideals of finite
co-length. In rings such as $k[s,t]$, due to a theorem of Serre,
irreducible ideals are complete intersections, a fact that leads to
Sylvester forms of low degree.

Turning to the equations of almost  complete intersections,
 we derive several Sylvester forms over a
polynomial ring $R=k[s,t]$, package them into  ideals and examine
the incident homological properties of these ideals and the
associated algebras.
 It is a
computer-assisted approach whose  role is to  produce a set of
syzygies that afford hand computation: the required equations
themselves are {\em not} generated by computation. Concretely, we
model a generic class of ideals cases to define `super-generic'
ideals $L$ in rings with several new variables
\[ L = (f,g, h_1, \ldots, h_m)\subset A.\]
Using {\it Macaulay2} (\cite{Macaulay2}), we obtain the free
resolutions of $L$. In degrees $\leq 5$, the resolution has length
$\leq 3$ ($2$ when degree $=4$)
\[ 0 \rar F_3 \stackrel{d_3}{\lar}
 F_2 \stackrel{d_2}{\lar} F_1 \lar F_0 \lar L \rar 0.
\] It has the property that after specialization the ideals of
maximal minors of $d_3$ and $d_2$ have codimension $5$ and $\geq 4$,
respectively. Standard arguments of the theory of free resolutions
will suffice to show that the specialization of $L$ is a prime
ideal.

For ideals in $R=k[s,t]$ generated by forms of  degrees $\leq 5$,
the method succeeds in describing the full set of equations. In
higher degree, in cases of special interest, it predicts the precise
form of the elimination equation.

For a technical reason--due to the character of irreducible
ideals--the method is limited to dimension  two. Nevertheless, it is
supple enough to apply to non-homogeneous ideals. This may be
exploited elsewhere, along with the treatment of ideals with larger
numbers of generators in a two-dimensional ring.

\section{Preliminaries on symmetric and Rees algebras} \label{Prelim}

We will introduce some basic material of Rees algebras (\cite{BH},
 \cite{HSV1}, \cite{alt}). Since most of the
questions we will consider have a local character, we pick local
rings as our setting. Whenever required, the transition to graded
rings will be direct.

\medskip

Throughout we will consider
 a Noetherian local ring $(R, \mathfrak{m})$
and $I$ an $\mathfrak{m}$-primary ideal (or a graded algebra over a
field $k$, $R= \sum_{n\geq 0}R_n= R_0[R_1]$, $R_0=k$, and $I$ a
homogeneous ideal of finite colength $\lambda(R/I)< \infty$).

We assume that $I$ admits a minimal reduction $J$ generated by
$n=\dim R$ elements. This is always possible  when $k$ is infinite.
The terminology means that for some integer $r$, $I^{r+1}= JI^r$.
This condition in turn means that the inclusion of Rees algebras
$R[Jt]\subset R[It]$ is an integral birational extension (birational
in the sense that the two algebras have the same total ring of
fractions). The smallest such integer, $r_J(I)$, is called the {\it
reduction number } of $I$ relative to $J$; the infimum of these
numbers over all minimal reductions of $I$ is the (absolute)
reduction number $r(I)$ of $I$.

For any ideal, not necessarily $\mathfrak{m}$-primary, the
 {\em special fiber} of $R[It]$ -- or of $I$ by abuse of terminology -- is
the algebra $\mathcal{F}(I)= R[It]\otimes_R(R/\mathfrak{m})$. The
dimension of $\mathcal{F}(I)$ is called the {\it analytic spread} of
$I$, and denoted $\ell(I)$. When $I$ is $\mathfrak{m}$-primary,
$\ell(I)=\dim R $. A minimal reduction $J$ is generated by $\ell(I)$
elements, and $\mathcal{F}(J)$ is a Noether normalization of
$\mathcal{F}(I)$.

\subsection*{Hilbert polynomials}

 The Hilbert polynomial of $I$ by ($m\gg 0$) is the function
 (\cite{BH}):

\[ \lambda(R/I^m) = e_0(I){{m+n-1}\choose{n}} - e_1(I)
{{m+n-2}\choose{n-1}} + \textrm{lower terms}.\] $e_0(I)$ is the
multiplicity of the ideal $I$. If $R$ is Cohen-Macaulay, $e_0(I) =
\lambda(R/J)$, where $J$ is a minimal reduction of $I$ (generated by
a regular sequence). For such rings, $e_1(I)\geq 0$.

\medskip

For instance, if $R= k[x_1, \ldots, x_n]$,  $\mathfrak{m}= (x_1,
\ldots, x_n)$ and
 $I=\mathfrak{m}^d$,
\begin{eqnarray}\nonumber
\lambda(R/I^m) &=&
\lambda(R/\mathfrak{m}^{md})={{md+n-1}\choose{n}}\\ \nonumber
 &=& d^n{{m+n-1}\choose{n}}- e_1(I)
{{m+n-2}\choose{n-1}}+ \textrm{lower terms}
\end{eqnarray}
where $e_1(I)= {\frac{n-1}{2}}(d^n-d^{n-1})$.

\medskip

Both coefficients will be the focus of our interest soon.

\subsection*{Cohen-Macaulay Rees algebras}

There is broad array of criteria expressing the Cohen-Macaulayness
of Rees algebra (see \cite{AHT}, \cite{JK95}, \cite{SUV2},
\cite[Chapter 3]{icbook}). Our needs will be filled by single
criterion whose proof is fairly straightforward. We briefly review
its related contents.

Let $(R, \mathfrak{m})$ be a Cohen-Macaulay local ring of dimension
$\geq 1$, and let $I$ be an $\mathfrak{m}$-primary ideal with a
minimal reduction $J$. The Rees algebra $R[Jt]$ is Cohen-Macaulay
and serves as an anchor to derive many properties of $R[It]$. Here
is one that we shall make use of.

\medskip

Define the {\em Sally module} $S_J(I)$ of $I$ relative to $J$ to be
the cokernel of the natural inclusion of finite $R[Jt]$-modules $I\,
R[Jt]\subset I\, R[It]$. Thus,
\[ S_J(I) = \sum_{t\geq 2} I^t/IJ^{t-1}.\]
 It has a
Hilbert function, unlike the algebra $R[It]$, that gives information
about the Hilbert function of $I$ (see \cite[Chapter 2]{alt}). The
module on the left, $I\cdot R[Jt]$, is a Cohen-Macaulay
$R[Jt]$-module of depth $\dim R+1$. The Cohen-Macaulayness of
$I\cdot R[It]$ is directly related to that of $R[It]$. These
considerations lead to the criterion:

\begin{Theorem} \label{cmtest} If $\dim R\geq 2$ and
  the reduction number of $I$ is
$\leq 1$, that is $I^2=JI$, then $R[It]$ is Cohen-Macaulay. The
converse holds if $\dim R=2$.
\end{Theorem}

\subsection*{Symmetric algebras}

Throughout $R$ is a Cohen-Macaulay ring and $I$ is an almost
complete intersection. The symmetric algebra $\Sym(I)$ will be
denoted by $\Symi$. Hopefully there will be no confusion between
$\Symi$  and the rings of polynomials $S=R[T_1, \ldots, T_n]$ that
we use to give a presentation of either $\Rees$ or $\Symi$.

\medskip

What keeps symmetric algebras of almost complete intersections
fairly under control is the following:

\begin{Proposition} Let $(R, \mathfrak{m})$ be a Cohen-Macaulay local ring. If $I$ is
an almost complete intersection and $\depth R/I\geq \dim R/I-1$,
then $\Symi$ is Cohen-Macaulay. In particular, if $I$ is
$\mathfrak{m}$-primary then $\Symi$ is Cohen-Macaulay.
\end{Proposition}

\demo The general assertion follows from \cite[Proposition
10.3]{HSV1}; see also \cite{Rossi}. \QED

\bigskip

Let $R$ be a Noetherian  ring and let $I$ be an $R$-ideal with a
free presentation
\[ R^m \stackrel{\varphi}{\lar} R^n \lar I \rar 0.\]
We assume that $I$ has a regular element. If $S=R[T_1, \ldots,
T_n]$, the symmetric algebra $\Symi$ of $I$ is defined by the ideal
$M_1\subset S$ of $1$-forms,
\[ M_1= I_1([T_1 , \ldots, T_n] \cdot \varphi).
\]
The ideal of definition of the Rees algebra $\Rees$ of $I$ is the
ideal $M\subset S$ obtained by elimination
\[ M= \bigcup_t (M_1: x^t)=M_1:
x^{\infty},
\] where $x$ is a regular element of $I$.

\subsection*{Sylvester forms}

To get additional elements of $M$, evading the above calculation, we
make use of general Sylvester forms. Recall how
 these are obtained. Let
$\mathbf{f}= \{f_1, \ldots, f_n\} $ be a set of polynomials in
 $B=R[x_1, \ldots, x_r]$ and let $\mathbf{a}= \{a_1, \ldots,
a_n\}\subset R $. If  $f_i \in (\mathbf{a})B$ for all $i$, we can
write
\[ \mathbf{f}= [f_1 \cdots f_n] = [a_1 \cdots a_n] \cdot A =
\mathbf{a}\cdot A,\] where $A $ is a $n\times n$ matrix with entries
in $B$. By an abuse of terminology, we refer to $\det(A)$ as a {\em
Sylvester form} of $\mathbf{f}$ relative to $\mathbf{a}$, in
notation
\[ \det(\mathbf{f})_{(\mathbf{a})}=\det(A).\]
It is not difficult to show that $ \det(\mathbf{f})_{(\mathbf{a})}$
is well-defined mod $(\mathbf{f})$. The classical Sylvester forms
are defined relative to sets of monomials (see \cite{Cox}). We will
make use of them in Section~\ref{dimtwo}. The structure of the
matrix $A$ may give rise to finer constructions (lower order
Pfaffians, for example) in exceptional cases (see \cite{dual}).

In our approach, the $f_i$ are elements of $M_1$,  or were obtained
in a previous calculation, and the ideal $(\mathbf{a})$ is derived
from the matrix of syzygies $\varphi$.

\section{Algebraic invariants in rational parametrizations} \label{Ratparam}

 Let  $f_1, \ldots, f_{n+1}\in R=k[x_1, \ldots, x_n]$ be forms of
the same degree. They define a rational map
\[\Psi: \mathbb{P}^{n-1} \dasharrow \mathbb{P}^n \]
\[p \rightarrow (f_1(p): f_2(p): \cdots : f_{n+1}(p)).\]
Rational maps are defined more generally with any number $m$ of
forms of the same degree, but in this work we only deal with the
case where $m=n+1$.

There are two basic ingredients to the algebraic side of rational
map theory: the ideal theoretic  and the algebra aspects, both
relevant for the nature of $\Psi$. First the ideal $I=(f_1, \ldots,
f_{n+1})\subset R$, which in this context is called the {\em base
ideal} of the rational map. Then there is the $k$-subalgebra $k[f_1,
\ldots, f_{n+1}]\subset R$, which is homogeneous, hence a standard
$k$-algebra up to degree renormalization. As such it gives the
homogeneous coordinate ring of the (closed) image of $\Psi$. Finding
the irreducible defining equation of the image is known as {\em
elimination} or {\em implicitization}.

We refer to \cite{ram2} and \cite{bir2003} (also \cite{dual} for an
even earlier overview) for the interplay between the ideal and the
algebra, as well as its geometric consequences. In particular, the
Rees algebra $\Rees=R[It]$ plays a fundamental role in the theory. A
pleasant side of it is that, since $I$ is generated by forms of the
same degree, one has $\Rees\otimes_R k\simeq k[f_1t, \ldots,
f_{n+1}t]\subset \Rees$, which retro-explains  the (closed) image of
$\mathbb{P}^{n-1}$ by $\Psi$ as the image of the projection to
$\mathbb{P}^{n}$ of the graph of $\Psi$. In particular, the fiber
cone is reduced and irreducible.



\subsection{Elimination degrees and birationality}

Although a rational map $\pp^{n-1}\dasharrow \pp^n$ has a unique set
of defining forms $f_1,\ldots,f_{n+1}$ of the same degree and unit
gcd, two such maps may look ``nearly'' the same if they happen to be
composite with a birational map of the target $\pp^n$ - a so-called
Cremona transformation. If this is the case the two maps have the
same degree, in particular the final elimination degrees are the
same.

However, it may still be the case that the two maps are composite
with a rational map of the target which is not birational, so that
their degrees as maps do not coincide, yet the degrees of the
respective images are the same. In such an event, one would like to
pick among all such maps one with smallest possible degree. This
leads us to he notion of improper and proper rational
parametrizations.

\begin{Definition}\label{def_of_proper}\rm
Let $\Psi=(f_1:\cdots :f_{n+1}): \pp^{n-1}\dasharrow \pp^n$ be a
rational map, where $\gcd(f_1, \ldots, f_{n+1})=1$. We will say that
$\Psi$ (or the parametrization defined by $f_1, \ldots, f_{n+1}$) is
{\em improper} if there exists a rational map
$$\Psi'=(f_1':\cdots :f_{n+1}'): \pp^{n-1}\dasharrow \pp^n,$$
with $\gcd(f_1', \ldots, f_{n+1}')=1$, such that:
\begin{enumerate}
\item There is an inclusion of $k$-algebras $k[f_1, \ldots, f_{n+1}]
\subset k[f_1', \ldots, f_{n+1}']$;
\item There is an isomorphism of $k$-algebras $k[f_1, \ldots, f_{n+1}]
\simeq k[f_1', \ldots, f_{n+1}']$;
\item $\deg \Psi' < \deg \Psi$.
\end{enumerate}
\end{Definition}
We note that if $\Psi$ is improper and $\Psi'$ is as above then the
rational map
$$(P_1:\cdots:P_{n+1}):\pp^n\dasharrow \pp^n$$
 is not birational, where $f_j=P_j(f_1', \ldots, f_{n+1}')$, for $1\leq
j\leq n+1$. Of course, the transition forms $P_j=P_j(y_1,\ldots,
y_{n+1})$ are not uniquely defined.

\begin{Example}\rm The parametrization given by
$f_1=x_1^4,f_2=x_1^2x_2^2, f_3=x_2^4$ is improper since it factors
through the parametrization $f'_1=x_1^2,f'_2=x_1x_2, f'_3=x_2^2$
through either one of the rational maps $(y_1:y_2:y_3)\mapsto
(y_1^2:y_2^2:y_3^2)$ or $(y_1:y_2:y_3)\mapsto (y_1^2:y_1t_3:y_3^2)$
neither of which is birational. Moreover, the forms $x_1^2,x_1x_2,
x_2^2$ define a birational map onto its image.
\end{Example}

We say that a rational map $\Psi=(f_1:\cdots :f_{n+1}):
\pp^{n-1}\dasharrow \pp^n$ is {\em proper} if it is not improper.
The need for considering proper rational maps will become apparent
in the context. It is also a basic assumption in elimination theory
when one is looking for the elimination degrees (see \cite{Cox}).

Clearly, if $\Psi$ is birational onto its image then it is proper.
The converse does not hold and one seeks for precise conditions
under which $\Psi$ is birational onto its image. This is the object
of the following parts of this subsection.

\medskip

When the ideal $I=(f_1, \ldots, f_{n+1})$ has finite co-length --
that is, $I$ is  $(x_1, \ldots, x_n)$-primary -- it is natural to
consider
 another mapping, namely, the corresponding embedding of the
Rees algebra $\Rees=R[It]$ into its integral closure
$\widetilde{\Rees}$. We will  explore the attached Hilbert functions
into the determinations of various degrees, including the
elimination degree of the mapping.

Thus, assume that $I$ has finite co-length. Then we may assume ($k$
is infinite) that $f_1, \ldots, f_n$ is a regular sequence, hence
the multiplicity of $J=(f_1, \ldots, f_n)$ is $d^n$, the same as the
multiplicity of $\mathfrak{m}^d$. This implies that $J$ is a minimal
reduction of $I$ and of $\mathfrak{m}^d$. We will set up a
comparison between $\Rees$ and $\Rees'=R[\mathfrak{m}^d]$, where
$\mathfrak{m}=(x_1, \ldots, x_n)$, through two relevant exact
sequences:
\begin{eqnarray} \label{es1}
0 \rar \Rees \lar \Rees' \lar D \rar 0,
\end{eqnarray}
and its reduction mod $\mathfrak{m}$
\begin{eqnarray} \label{es2}
\bar{\Rees} \lar \bar{\Rees'} \lar \bar{D} \rar 0.
\end{eqnarray}
$\F=\bar{\Rees}$ is the {\em special fiber} of $\Rees$ (or, of $I$),
and since $I$ is generated by forms of the same degree, one has
$\F\simeq k[f_1, \ldots, f_{n+1}]$ as graded $k$-algebras. By the
same token,  $\F'=\bar{\Rees'}\simeq  k[\mathfrak{m}^d]$ -- the
$d$-th Veronese subring of $R$. In particular, since $\dim \F=\dim
\F'$, the leftmost map in the exact sequence (\ref{es2}) is
injective. Also $D$ is annihilated by a power of $\mathfrak{m}$,
hence $\dim D= \dim \bar{D}$.

These are the degrees (multiplicities)
 $\deg(\F)$
and $\deg(\F')$ of the special fibers. Since $\F'$ is an integral
extension of $\F$, one has
\begin{eqnarray} \label{es3}
\deg(\F') = \deg(\F) [\F':\F],
\end{eqnarray}
where $[\F':\F]=\dim _K(\F'\otimes _\F K)$, where $K$ denotes the
fraction field of $\F$ (see, e.g., \cite[Proposition 6.1 (b) and
Theorem 6.6]{ram2} for more general formulas). Since $\F'$ is
besides integrally closed, the latter is also the field extension
degree $[\,k({\mathfrak m}^d)\,:\,K\,]$. Note that $[\F':\F]=1$
means that the extension $\F\subset \F'$ is birational
(equivalently, the rational map $\Psi$ maps ${\mathbb P}^{n-1}$
birationally onto its image). As above, set $L={\mathfrak m}^d$. We
next characterize birationality in terms of both the coefficient
$e_1$ and the dimension of the ${\mathcal R}$-module $D$.

\begin{Proposition}\label{birationality}
The following conditions are equivalent:
\begin{enumerate}
\item[{\rm (i)}] $[\F':\F]=1$, that is $\Psi$ is birational onto
its image;
\item[{\rm (ii)}] $\deg (\F)=d^{n-1}$;
\item[{\rm (iii)}] $\dim \bar{D}\leq n-1$;
\item[{\rm (iv)}] $\dim D\leq n-1$
\item[{\rm (v)}] $e_1(L)=e_1(I)$.
\end{enumerate}
\end{Proposition}
\demo (i) $\Longleftrightarrow$ (ii) This is clear from (\ref{es3})
since $\deg(\F')=d^{n-1}$.

(i) $\Longleftrightarrow$ (iii) Since $\ell(I)=n$ and $\F\subset
\F'$ is integral, then $\F\subset \F'$ is a birational extension if
and only if its conductor $\F:_\F\F'$ is nonzero, equivalently, if
and only if $\dim \bar{D}\leq n-1$.


(iv) $\Longleftrightarrow$ (iii) Clearly, $\dim D\leq n$ and in the
case of equality its multiplicity is $e_1(L)-e_1(I)>0$. Therefore,
the equivalence of the two statements follows suit. \QED

\bigskip

There is some advantage in examining $\bar{D}$ since $\F$ is a
hypersurface ring,
\[ \F=k[T_1,\ldots, T_{n+1}]/(f)= R[T_1, \ldots, T_{n+1}]/(x_1,\ldots,
x_n,f)\] a complete intersection. Since $\F'$ is also
Cohen-Macaulay, with  a well-known presentation, it affords an
understanding of $\bar{D}$, and sometimes, of $D$.

\subsection{Calculation of $e_1(I)$ of the base ideal of a rational
map}

One objective here is to apply some general formulas for the Chern
number $e_1(I)$ of an ideal $I$ to the case of the base ideal of a
rational map with source $\pp^1={\rm Proj}(k[x_1,x_2])$.




\medskip

Here is a method put together from scattered facts in the literature
of Rees algebras (see  \cite[Chapter 2]{icbook}).

\begin{Proposition}\label{obs1} {\rm
Let $(R,\mathfrak{m})$ be a Cohen-Macaulay local ring of dimension
$d$, let $I$ be an $\mathfrak{m}$-primary ideal with a minimal
reduction $J=(a_1, \ldots, a_{d})$. Set $R'=R/(a_1, \ldots,
a_{d-1})$, $I' = IR'$. Then
\begin{itemize}
\item[{\rm (i)}] $e_0(I)=e_0(I')=\lambda(R/J)$, $e_1(I)=e_1(I')$
\item[{\rm (ii)}] $\red(I')< \deg R'\leq e_0(I)$; in particular,
for $n\geq r=\red(I')$, one has $I'^{\,n+1}=a_{d}I'^{\,n}$
\item[{\rm (iii)}] $\lambda(R'/I'^{\,r+1})= \lambda(R'/I'^{\,r})+
\lambda(I'^{\,r}/a_{d}I'^{\,r})= e_0(I)(r+1)-e_1(I)$
\item[{\rm (iv)}] $e_1(I)= -\lambda(R'/I'^{\,r})+ e_0(I)r$
\end{itemize}
}\end{Proposition}

It would be desirable to develop a direct  method suitable  for the
ideal $I=(a,b,c)$ generated by forms of $R=k[s,t]$, of degree $n$.
We may assume that $a,b$ for a regular sequence (i.e.
$\gcd(a,b)=1$).
  We
already know that $e_0(I)= n^2$. For regular rings, one knows
(\cite{ni1}) that $e_1(I)\leq {\frac{d-1}{2}}\,e_0(I)$, $d=\dim R$.
Nevertheless the steps above already lead to an efficient
calculation for two reasons: the multiplicity $e_0(I)$ is known at
the outset  and  it does not really involve the powers of $I$. Forms
of degree up to 10 are handled well by {\em Maucalay2}
(\cite{Macaulay2}).

\section{Sylvester forms in dimension two}\label{dimtwo}

We establish the basic notation to be used throughout. $R= k[s,t]$
is a polynomial ring over the infinite field $k$, and $I\subset R=
k[s,t]$ is a codimension $2$  ideal generated by $3$
 forms of the same degree $n+1$, with free graded resolution
\[
0  \longrightarrow R(-n-1-\mu) \oplus R(2(-n-1)+\mu)
\stackrel{\varphi}{\longrightarrow} R^3(-n-1) \longrightarrow I
\longrightarrow 0, \quad \varphi=
\left[\begin{array}{ccc}\alpha_1 & \beta_1 & \gamma_1 \\
\alpha_2 & \beta_2 & \gamma_2
\end{array}\right]^t.
\]

Then the symmetric algebra of $I$ is $\Symi\simeq R[T_1, T_2,
T_3]/(f,g)$ with
\begin{eqnarray*} \label{genf}
f &=& \alpha_1T_1+ \beta_1T_2+ \gamma_1T_3 \\ \label{geng} g &=&
\alpha_2T_1+ \beta_2T_2+ \gamma_2T_3.
\end{eqnarray*}
Starting out from  these $2$ forms,
  the defining equations
of $\Symi$, following \cite{Cox}, we obtain by elimination higher
degrees forms in the defining ideal of $\Rees(I)$. It will  make use
of a computer-assisted methodology to show that these
algorithmically specified sets generate the  ideal of definition $M$
of $\Rees(I)$ in several cases of
 interest--in particular answering some questions raised
 \cite{Cox}. More precisely, the so-called ideal of moving forms $M$
 is given when $I$ is generated by forms of degree at most $5$. In
 arbitrary degree, the algorithm  will provide the elimination equation in
 significant cases.

\subsection{Basic Sylvester forms in dimension $2$}

Let $R=k[s,t]$ and let $F,G\in B= R[s,t,T_1,T_2,T_3]$. If $F,G\in
(u,v)B$, for some ideal $(u,v)\subset R $, the form derived from
\[
\left[ \begin{array}{r}
f \\g \\
\end{array} \right]=
\left[ \begin{array}{rr}
a & b \\c & d \\
\end{array} \right]
\left[ \begin{array}{r}
u \\v \\
\end{array} \right],\]
\[ h=ad-bc = \det(F,G)_{(u,v)}, \]
will be called a basic Sylvester form.

To explain their naturalness, even for ideals $I$ not necessarily
generated by forms, we give an approach to irreducible decomposition
of certain ideals.

\begin{Theorem} Let $(R,\mathfrak{m})$ be a Gorenstein local ring and let
$I$ be an $\mathfrak{m}$--primary ideal. Let $J\subset I$ be an
ideal generated by a system of parameters and let $E = (J:I)/J$ be
the canonical module of $R/I$. If $E = (e_1, \ldots, e_r)$, $e_i\neq
0$, and $I_i= \ann(e_i)$, then $I_i$ is an irreducible ideal and
\[ I = \bigcap_{i=1}^r I_i.\]
\end{Theorem}

The statement and its proof will apply to ideals of rings of
polynomials over a field.

\bigskip

\demo The module $E$ is the injective envelope of $R/I$, and
therefore it is a faithful $R/I$--module (see \cite[Section 3.2]{BH}
for relevant notions). For each $e_i$, $Re_1$ is a nonzero submodule
of $E $ whose socle is contained in the socle of $E$ (which is
isomorphic to $R/\mathfrak{m}$) and therefore its annihilator $I_i$
(as an $R$-ideal) is irreducible. Since the intersection of the
$I_i$ is the annihilator of $E$, the asserted equality follows. \QED

\begin{Corollary} Let $(R, \mathfrak{m})$ be a regular local ring of
dimension two and let $I$ be an $\mathfrak{m}$--primary ideal with a
free resolution
\[ 0 \rar R^{n-1} \stackrel{\varphi}{\lar} R^n \lar I \rar 0,\]
\[ \varphi = \left[ \begin{array}{lcr}
& & \\
& \varphi' & \\
& & \\
\hline
a_{n-1,1} & \cdots & a_{n-1,n-1} \\
a_{n,1} & \cdots & a_{n,n-1} \\
\end{array} \right],
\]
 and suppose that the last two maximal minors $\Delta_{n-1}, \Delta_n$ of $\varphi$
 form a regular sequence. If $e_1, \ldots, e_{n-1}$ are as above,
 then
\[  (\Delta_{n-1}, \Delta_n): I = I_{n-2}(\xi') = (e_1, \ldots, e_{n-1})\]
and each ideal $(\Delta_{n-1}, \Delta_n):e_i$ is a complete
intersection of codimension $2$.
\end{Corollary}

\demo The assertion that the irreducible $I_i$ is a complete
intersection is a result of Serre, valid for all two-dimensional
regular  rings whose projective modules are free. \QED

%
%
%
%
%
%
%
%
%
%

\begin{Remark}{\rm
In our applications, $I= C(f,g)$,
 the content ideal of $f,g$.
In some of these cases, $C(f,g) = (s,t)^n$, for some $n$, an ideal
which admits the irreducible decomposition
\[ (s,t)^n = \bigcap_{i=1}^n (s^i, t^{n+1-i}).\]
One can then
 process $f,g$ through all the pairs $\{s^i, t^{n-i+1}\}$,
and collect  the determinants for the next round of elimination.  As
in the classical Sylvester forms, the inclusion $C(f,g)\subset
(s,t)^n$ may be used anyway to start the process,  although without
the measure of control of degrees afforded  by the  equality of
ideals. }\end{Remark}

\subsection{Cohen-Macaulay algebras}

We pointed out in Theorem~\ref{cmtest} that the
 basic control of Cohen-Macaulayness of a Rees algebra
of an ideal $I\subset  k[s,t]$ is that its reduction number be at
most $1$. We next give a mean of checking this property directly off
a free presentation of $I$.

\begin{Theorem}\label{macaulay_case} Let $I\subset R$ be an ideal  of codimension $2$, minimally
generated by $3$ forms of the same degree. Let
\[ \varphi= \left[ \begin{array}{rr}
\alpha_1 & \alpha_2 \\
\beta_1 & \beta_2 \\
\gamma_1 & \gamma_2\\
\end{array} \right] \]
be the Hilbert-Burch presentation matrix of $I$. Then $\Rees$ is
Cohen-Macaulay if and only if the equalities of ideals of $R$ hold
\[ (\alpha_1, \beta_1, \gamma_1)=(\alpha_2, \beta_2, \gamma_2)=(u,v),
\] where $u,v$ are forms.
\end{Theorem}

\demo Consider  the presentation
\[ 0 \rar \mathcal{L} \lar \Symi=R[T_1, T_2, T_3]/(f,g) \lar \Rees \rar 0,\]
where $f,g$ are the $1$-forms
\[
\left[
\begin{array}{r}
f \\ g
\end{array} \right]
= \left[
\begin{array}{rrr}
T_1 & T_2 & T_3 \\
\end{array} \right] \cdot \varphi.
\]

If $\Rees$ is Cohen-Macaulay, the reduction number of $I$ is $1$ by
Theorem~\ref{cmtest}, so there must be a nonzero quadratic form $h$
with coefficients in $k$ in  the presentation ideal $M$ of $\Rees$.
In addition to $h$, this ideal contains $f,g$,  hence in order to
produce such terms its Hilbert-Burch matrix must be of the form
\[
\left[
\begin{array}{ll}
u & v \\
 p_1 & p_2 \\
q_1 & q_2 \\
\end{array} \right]
\] where $u,v$ are forms of $k[s,t]$, and the other entries are
$1$-forms of $k[T_1,T_2,T_3]$. Since $p_1,p_2$ are $q_1,q_2$ are
pairs of linearly independent $1$-forms, the assertion about the
ideals defined by the columns of $\varphi$ follow.

\subsection{Base ideals generated in degree $4$}

This is the case treated by D. Cox in his Luminy lecture
(\cite{Cox}). We accordingly change the notation to $R=k[s,t]$, $I=
(f_1, f_2, f_3)$, forms of degree $4$. The field $k$ is infinite,
and we further assume that $f_1, f_2$ form a regular sequence so
that $J=(f_1, f_2)$ is a reduction of $I$ and of $(s,t)^4$. Let
\begin{eqnarray}\label{es4}
0 \rar R(-4-\mu)\oplus R(-8+\mu)  \stackrel{\varphi}{\lar} R^3(-4)
\lar R\lar R/I\rar 0, \quad
 \varphi= \left[ \begin{array}{rr}
\alpha_1 & \alpha_2 \\
\beta_1 & \beta_2 \\
\gamma_1 & \gamma_2\\
\end{array} \right]
\end{eqnarray}
 be the Hilbert-Burch presentation of $I$. We
obtain the equations of $f_1, f_2, f_3$ from this matrix.

\medskip

Note that $\mu$ is the degree of the first column of $\varphi$,
$4-\mu$ the other degree. Let us first consider (as in \cite{Cox})
the case $\mu=2$.

\subsection*{Balanced case}

We shall now give a computer-assisted treatment of the {\em
balanced} case, that is when the resolution (\ref{es4}) of the ideal
$I$ has $\mu=2$  and the content ideal of the syzygies is $(s,t)^2$.
Since $k$ is infinite, it is easy to show that there is a change of
variables, $T_1, T_2, T_3 \rar x,y,z$, so that $(s^2, st, t^2)$ is a
syzygy of $I$. The forms $f,g$ that define the symmetric algebra of
$I$ can then be written
\[ [f \quad g] =
 [s^2 \quad st \quad t^2] \left[\begin{array}{rr} x & u \\
y & v \\
z & w \\
\end{array} \right],
\]
where $u,v,w$ are linear forms in $x,y,z$. Finally, we will assume
that
 the ideal
 $I_2\left(\left[ \begin{array}{rrr}
x & y & z \\
u & v & w \\
\end{array} \right]\right)$ has codimension two.
Note that this is a generic condition.


\medskip

We introduce now the {\em equations} of $I$.

\bigskip

 $\bullet$ Linear equations $f$ and $g$:

\begin{eqnarray*}
[f \quad g ] &=& [x \; y \; z]\; \varphi = [x \; y \; z] \left[
\begin{array}{rr}
\alpha_1 & \alpha_2 \\
\beta_1 & \beta_2 \\
\gamma_1 & \gamma_2\\
\end{array} \right] \\
&=& [s^2 \quad st \quad t^2]  \left[\begin{array}{rr} x & u \\
y & v \\
z & w \\
\end{array} \right],
\end{eqnarray*}
where $u,v,w$ are linear forms in $x,y,z$.

\bigskip

{ $\bullet$ Biforms $h_1$ and $h_2$:}

\bigskip

Write $\Gamma_1 $ and $\Gamma_2$ such that
\[
[f \quad g ]=[x \; y \; z]\; \varphi =[\;s \quad t^2\;]\;
 \Gamma_1 = [\;s^2 \quad t\;]\; \Gamma_2.
\]

Then $h_1=\det \Gamma_1$ and $h_2=\det \Gamma_2$.

\bigskip
{ $\bullet$ Implicit equation $F=\det \Theta$,
 where $[ h_1 \quad h_2] = [s \quad t]\; \Theta$.}

\bigskip

Using generic entries for $\varphi$, in place of the true $k$-linear
forms in old variables $x,y,z$, we consider the ideal of
$k[s,t,x,y,z,u,v,w]$ defined by

\begin{eqnarray*}  f &=& s^2x+sty+t^2z \\
 g &=&  s^2u+stv+t^2w\\
h_1&=& - syu - tzu + sxv + txw\\
h_2 &=& - szu - tzv + sxw + tyw\\
F &=& -z^2u^2 + yzuv-xzv^2-y^2uw+2xzuw+xyvw-x^2w^2
\end{eqnarray*}

\begin{Proposition} \label{case24}
If $I_2\left(\left[ \begin{array}{rrr}
x & y & z \\
u & v & w \\
\end{array} \right]\right)
$ specializes to a codimension two ideal of $k[x,y,z]$, then
 $ L = (f,g,h_1,h_2,F)\subset A = R[x,y,z,u,v,w]$
 specializes
to the defining ideal of $\Rees$.
\end{Proposition}


\medskip

\demo {\em Macaulay2} (\cite{Macaulay2}) gives a resolution
\[ 0 \rar A \stackrel{d_2}{\lar} A^5 \lar A^5 \lar L \rar 0\]
where
\[ d_2 = \left[
\begin{array}{c}
zv-yw \\
zu-xw \\
-yu+xv \\
\kern-2pt-t \\
s\\
\end{array} \right].
\]
The assumption on
 $I_2\left(\left[ \begin{array}{rrr}
x & y & z \\
u & v & w \\
\end{array} \right]\right)$ says that the entries of $d_2$ generate
an ideal of codimension four and thus implies that the
specialization $LS$ has projective dimension two and that it is
unmixed. Since $LS \not \subset (s,t)S$,
 there is an element  $q\in (s,t)R$
 that is regular modulo $S/LS$. If
\[ LS = Q_1 \cap \cdots \cap Q_r\]
is the primary decomposition of $LS$, the
 localization
 $LS_q$ has the corresponding decomposition since $q$ is not contained
 in any of the $\sqrt{Q_i}$. But now $\Sym_q=\Rees_q$, so
 $LS_q=(f,g)_u$, as $I_q=R_q$.
\QED

\subsection*{Non-balanced case}


We shall now give a similar computer-assisted treatment of the
non-balanced case, that is when the resolution (\ref{es4}) of the
ideal $I$ has $\mu=3$. This implies that the content ideal of the
syzygies is $(s,t)$. Let us first indicate how the proposed
algorithm would behave.

\begin{itemize}

\item[{$\bullet$}] Write the forms $f,g$ as
\begin{eqnarray*}
f & = & a s + b t \\
g & = & c s + d t, \\
\end{eqnarray*}
where
\[
 \left[ \begin{array}{r}
c \\[5pt]
d \\
\end{array}\right]
=
 \left[ \begin{array}{rrr}
x & y & z \\[5pt]
u & v & w \\
\end{array}\right]
 \left[ \begin{array}{l}
s^2 \\[5pt]
st \\[5pt]
t^2
\end{array}\right]
\]

\item[{$\bullet$}] The next form is the Jacobian of $f,g$ with
respect to $(s,t)$
\[ h_1 = \det(f,g)_{(s,t)}= ad - bc =
 - bxs^2 - byst - bzt^2 + aus^2  + avst + awt^2.
\]
\item[{$\bullet$}] The next two generators
\[ h_2 = \det(f,h_1)_{(s,t)} =
 b^2x s + b^2y t - abzt - abus - abvt + a^2w t
\]
and the elimination equation
\[ h_3 = \det(f,h_2)_{(s,t)} =
 - b^3 x + ab^2 y - a^2 bz + ab^2 u - a^2 bv + a^3 w.
\]


\end{itemize}

\begin{Proposition} \label{case41}
$ L = (f,g,h_1,h_2,h_3)\subset A = k[s,t, x,y,z,u,v,w]$
  specializes
to the defining ideal of $\Rees$.
\end{Proposition}

\demo {\em Macaulay2} (\cite{Macaulay2}) gives the following
resolution of $L$
\[ 0 \rar A^2 \stackrel{\varphi}{\lar} A^6 \stackrel{\psi}{\lar} A^5
\lar L \rar 0,\]
\begin{scriptsize}
\[
\varphi = \left[
\begin{array}{rr}
s & 0 \\
t & 0 \\
-b & s \\
a & t \\
0 & -b \\
0 & a \\
\end{array}
\right], \]
\end{scriptsize}
\begin{tiny}
\[ \psi =
\left[
\begin{array}{cccccc}
-b^2x+abu & -b^2y+abz+abv-a^2w & -bsx-bty+asu+atv & -btz+atw &
   -s^2x-sty-t^2z & -s^2u-stv-t^2w \\[5pt]
t         &       -s           &        0         &    0     &
      0           &         0      \\[5pt]
a         &        b           &        t         &    -s    &
      0           &         0      \\[5pt]
0         &        0           &        a         &     b    &
      t           &        -s      \\[5pt]
0         &        0           &        0         &     0    &
      a           &         b      \\
\end{array} \right]
\]
\end{tiny}

The ideal of $2\times 2$ minors of $\varphi$ has codimension $4$,
even after we specialize from $A$ to $S$ in the natural manner.
Since $LS$ has projective dimension two, it will be
 unmixed. As $LS\not \subset (s,t)$, there is an element  $u\in (s,t)R$
 that is regular modulo $S/LS$. If
\[ LS = Q_1 \cap \cdots \cap Q_r\]
is the primary decomposition of $LS$, the
 localization
 $LS_u$ has the corresponding decomposition since $u$ is not contained
 in any of the $\sqrt{Q_i}$. But now $\Sym_u=\Rees_u$, so
 $LS_u=(f,g)_u$, as $I_u=R_u$.
 \QED

\subsection{Degree $5$ and above}

It may be worthwhile to extend this to arbitrary degree, that is
assume that $I$ is defined by $3$ forms of degree $n+1$ (for
convenience in the notation to follow). We first consider the case
$\mu=1$.
 Using the procedure above, we
would obtain the  sequence of polynomials in $A= R[a,b,x_1, \ldots,
x_n, y_1,\ldots, y_n]$

\begin{itemize}

\item[{$\bullet$}] Write the forms $f,g$ as
\begin{eqnarray*}
f & = & a s + b t \\
g & = & c s + d t, \\
\end{eqnarray*}
where
\[
 \left[ \begin{array}{r}
c \\
d \\
\end{array}\right]
=
 \left[ \begin{array}{rrr}
x_1 & \cdots  & x_n \\
y_1 & \cdots  & y_n \\
\end{array}\right]
 \left[ \begin{array}{c}
s^{n-1} \\
s^{n-2}t \\
\vdots \\
st^{n-2}\\
t^{n-1}
\end{array}\right]
\]

\item[{$\bullet$}] The next form is the Jacobian of $f,g$ with
respect to $(s,t)$
\[ h_1 = \det(f,g)_{(s,t)}= ad - bc
\]
\item[{$\bullet$}] Successively we would set
\[ h_{i+1} = \det(f,h_i)_{(s,t)}, \quad 1 < n.
\]
\item[{$\bullet$}]  The polynomial
\[h_{n} = \det(f,h_{n-1})_{(s,t)}\]
is   the elimination equation.


\end{itemize}

\begin{Proposition} \label{case15} $ L = (f,g,h_1,\ldots, h_5)\subset A $
 specializes
to the defining ideal of $\Rees$.
\end{Proposition}

In {\em Macaulay2}, we checked the degrees $5$ and $6$ cases. In
both cases, the ideal $L$ (which has one more generator in degree
$6$) has a projective resolution of length $2$ and the ideal of
maximal minors of the last map has codimension four.

\begin{Conjecture} \label{mu1dimtwo}{\rm For arbitrary $n$,
$ L = (f,g,h_1,\ldots, h_n)\subset A $ has projective dimension two
and
  specializes
to the defining ideal of $\Rees$. }\end{Conjecture}

In degree $5$, the interesting case is when the Hilbert-Burch matrix
$\phi$ has degrees $2$ and $3$. Let us describe the proposed
generators. For simplicity, by a change of coordinates, we assume
that the coordinates of the degree $2$ column of $\varphi$ are
$s^2,st,t^2\,$

\begin{eqnarray*}
f & = & s^2x+sty+t^2z\\
g & = &
(s^3w_1+s^2tw_2+st^2w_3+t^3w_4)x+(s^3w_5+s^2tw_6+st^2w_7+t^3w_8)y\\
& + & (s^3w_9+s^2tw_{10}+st^2w_{11}+t^3w_{12})z
\end{eqnarray*}

Let
\[\begin{array}{lllll}
\left[ \begin{array}{c} f \\[5pt]
 g \end{array} \right] &=&
\left[ \begin{array}{lll} \kern3pt x &\quad\quad\quad y &\quad\;  z
\\[5pt]
sA &\quad  sB+tC &\quad tD
\end{array}
\right]
\left[ \begin{array}{c} s^2 \\
 st \\
  t^2 \end{array}
  \right] &=&
  \phi\left[ \begin{array}{c} s^2 \\
   st \\
    t^2 \end{array} \right] \\
    &&&&\\
&=& \left[ \begin{array}{ll}\kern15pt x &\kern17pt ys+zt  \\[5pt]
 sA+tB &\quad stC+t^2D
 \end{array} \right]
 \left[ \begin{array}{c} s^2 \\[5pt]
  t \end{array} \right]&=& B_1
  \left[ \begin{array}{c} s^2 \\[5pt]
   t \end{array} \right] \\&&&&\\
&=& \left[ \begin{array}{ll}\kern8pt xs+yt &\kern25pt z  \\[5pt]
 s^2A+stB &\quad
sC+tD \end{array} \right] \left[ \begin{array}{c} s \\[5pt]
 t^2
\end{array} \right] &=& B_2\left[ \begin{array}{c} s \\[5pt]
 t^2
\end{array} \right]
\end{array},
\]
where $A,B,C,D$ are $k$-linear forms in $x,y,z$.

\bigskip

\[\begin{array}{lll}
h_1 &= & \det (B_1) \\
&=& s^2(-yA)+st(xC-yB-zA)+t^2(xD-zB) \\
&=& s^2(-yA)+ t (xCs-yBs-zAs+xDt-zBt) \\
&=& s(-yAs+xCt-yBt-zAt )+t^2(xD-zB ), \\ &&\\
h_2 &=& \det (B_2) \\
&=& s^2(xC-zA)+st(xD+yC-zB)+t^2(yD) \\
&=& s^2(xC-zA)+t(xDs+yCs-zBs+yDt ) \\
&=& s(xCs-zAs+xDt+yCt-zBt )+t^2(yD ).
\end{array}
\]

\bigskip

\[\begin{array}{lllll}
\left[ \begin{array}{c} f \\[5pt]
 h_1 \end{array} \right] &=&
 \left[\begin{array}{ll} \kern12pt x &\kern62pt ys+zt \\
  -yA &\quad xCs-yBs-zAs+xDt-zBt  \end{array}\right]
  \left[ \begin{array}{c} s^2 \\[5pt]
   t \end{array} \right] &=& C_1
  \left[ \begin{array}{c} s^2 \\[5pt]
   t \end{array} \right] \\
    &&&&\\
&=&\left[\begin{array}{ll}\kern45pt xs+yt &\kern30pt z\\[5pt]
 -yAs+xCt-yBt-zAt &\quad
xD-zB   \end{array}\right]
\left[ \begin{array}{c} s \\[5pt]
 t^2 \end{array} \right]&=&C_2
 \left[ \begin{array}{c} s \\[5pt]
  t^2 \end{array} \right]
\end{array}
\]

\bigskip

\[\begin{array}{lllll}
\left[ \begin{array}{c} f \\[5pt]
 h_2 \end{array} \right] &=&
 \left[\begin{array}{ll} \kern18pt x &\kern62pt ys+zt \\[5pt]
  xC-zA  &\quad xDs+yCs-zBs+yDt \end{array}\right]
  \left[ \begin{array}{c} s^2 \\[5pt]
   t \end{array} \right] &=& C_3
   \left[ \begin{array}{c} s^2 \\[5pt]
    t \end{array} \right] \\ &&&&\\
&=&\left[\begin{array}{ll}\kern45pt xs+yt &\kern15pt z  \\[5pt]
 xCs-zAs+xDt+yCt-zBt &\quad yD   \end{array}\right]
 \left[ \begin{array}{c} s \\[5pt]
  t^2 \end{array} \right]&=&C_4
  \left[ \begin{array}{c} s \\[5pt]
   t^2 \end{array} \right]
\end{array}
\]

 \bigskip

 \[\begin{array}{lllll}
 c_1 &=& \det(C_1) &=& x^2(Cs+Dt)+xy(-Bs)+xz(-As-Bt)+yz(At)+y^2(As)\\ &&&&\\
 c_2 &=& \det(C_2) &=& x^2(Ds) + xy(Dt) + xz(-Bs-Ct)+yz(As)+z^2(At)\\&&&&\\
 c_3 &=& \det(C_3) &=& x^2(Ds) + xy(Dt) +xz(-Bs-Ct) + yz(As)+ z^2(At)\\ &&&&\\
 c_4 &=& \det(C_4) &=&  xy(Ds) + xz(-Cs-Dt) + yz(-Ct)+ z^2(As+Bt) + y^2(D)
 \end{array}
 \]

 \bigskip

  \[
 \left[ \begin{array}{c} f \\[5pt]
  h_1 \\
   h_2 \end{array} \right]=
 \left[ \begin{array}{lll}\kern15pt x &\kern38pt y &\kern30pt z \\[5pt]
 \kern5pt -yA &\quad xC-yB-zA &\quad xD-zB \\[5pt]
   xC-zA &\quad xD+yC-zB &\kern25pt yD  \end{array}\right]
   \left[ \begin{array}{c} s^2 \\[5pt]
    st \\[5pt]
     t^2 \end{array} \right]
 \]

 \bigskip

 Then $F= -x^3D^2+x^2yCD+xy^2(-BD) + x^2z(2BD-C^2)+xz^2(2AC-B^2)
 +xyz(BC-3AD)+y^2z(-AC)+yz^2(AB)+y^3(AD)+z^3(-A^2)$, an equation of degree $5$.
 In particular, the parametrization is birational.

\begin{Proposition} \label{case25}
 $L=(f,\, g,\, h_1,\, h_2,\, c_1,\, c_2,\, c_4,\, F)$
  specializes
to the defining ideal of $\Rees$.
\end{Proposition}

\demo Using {\it Macaulay2}, the ideal $L$ has a resolution:
\[
0 \longrightarrow S^1 \stackrel{d_3}{\longrightarrow} S^6
\stackrel{d_2}{\longrightarrow}  S^{12}
\stackrel{d_1}{\longrightarrow}   S^8 \longrightarrow L
\longrightarrow 0.
\]

\[
d_3= [-z \;\; y \;\; x \;\; -t \;\; s \;\; 0   ]^t
\]
\begin{scriptsize}
\[
d_2= \left[ \begin{array}{llllll} y & z&  0&  0 & 0&  \kern130pt 0 \\
  x & 0 & z&  0&  0 &  \kern130pt 0   \\
\kern-5pt -v & 0&  0&  z&  0&  \kern75pt x^2w_4-xzw_7+xyw_8+xzw_{12} \\
u&  0 & 0&  0 & z & -xzw_3+xyw_4+z^2w_6-yzw_7+y^2w_8-xzw_8-z^2w_{11}+yzw_{12} \\
 0 &  x& \kern-5pt -y& 0&  0& \kern130pt 0\\
 0 & \kern-5pt -v& 0& \kern-5pt -y& 0&  \kern60pt xzw_1-x^2w_3+yzw_5+z^2w_9-xzw_{11}\\
 0 & u & 0 & 0 &\kern-5pt -y & \kern70pt xzw_2-x^2w_4+z^2w_{10}-xzw_{12}  \\
 0 & 0&  u & 0& \kern-5pt -x & \kern60pt xzw_1+yzw_5-xzw_6+x^2w_8+z^2w_9\\
 0 & 0&  0&  u&  v& \kern130pt 0\\
 0 & 0&  v&  x&  0 & \kern15pt -xyw_1+x^2w_2-y^2w_5+xyw_6-x^2w_7-yzw_9+xzw_{10} \\
 0 & 0&  0&  0&  0& \kern125pt -t\\
 0 & 0&  0&  0&  0& \kern130pt s \end{array} \right]
\]
\end{scriptsize}

The ideals of maximal minors give $\codim I_1(d_3)=5$ and $\codim
I_{5}(d_2)= 4$ after specialization. As we have been arguing, this
suffices to show that the specialization is a prime ideal of
codimension two. \QED

\subsection*{Elimination forms in higher degree}

In degrees greater than $5$, the methods above are not very
suitable. However, in several cases they are still supple enough to
produce the elimination equation. We have already seen this when one
of the syzygies is of degree $1$. Let us describe two other cases.

\bigskip

$\bullet$ Degree $n=2p$, $f$ and $g$ both of degree $p$. We use the
decomposition
\[ (s,t)^p = \bigcap_{i=1}^p (s^i, t^{p+1-i}).\]
For each $1\leq i \leq p$, let
\[ h_i = \det(f,g)_{(s^i, t^{p+1-i})}.\]
These are quadratic polynomials with coefficients in $(s,t)^{p-1}$.
We set
\[ [h1, \cdots, h_p] = [s^{p-1}, \cdots, t^{p-1}]\cdot A,\]
where $A$ is a $p\times p$ matrix whose entries are $2$-forms in
$k[x,y,z]$. The Sylvester form of degree $n$, $F = \det(A)$, is the
required elimination equation.

\bigskip

$\bullet$ Degree $n=2p+1$, $f$  of degree $p$. We use the
decomposition
\[ (s,t)^p = \bigcap_{i=1}^p (s^i, t^{p+1-i}).\]
For each $1\leq i \leq p$, let
\[ h_i = \det(f,g)_{(s^i, t^{p+1-i})}.\]
These are quadratic polynomials with coefficients in $(s,t)^{p}$. We
set
\[ [f, h1, \cdots, h_p] = [s^{p}, \cdots, t^{p}]\cdot B,\]
where $A$ is a $(p+1)\times (p+1)$ matrix with one column whose
entries are linear forms and the remaining columns with entries
$2$-forms in $k[x,y,z]$. The Sylvester form $F = \det(B)$ is the
required elimination equation.

\end{document}